\magnification=\magstep0


\font\eu=eufm10 
\font\bigbf=cmbx12
\font\tenjap=msbm10\font\sevenjap=msbm7

\def\ts{\textstyle }

\def\bbC{{\mathchoice{\hbox{\tenjap C}}{\hbox{\tenjap C}}
{\hbox{\sevenjap C}}{\hbox{\sevenjap C}}}}
\def\bbO{{\mathchoice{\hbox{\tenjap O}}{\hbox{\tenjap O}}
{\hbox{\sevenjap O}}{\hbox{\sevenjap O}}}}
\def\bbR{ {\mathchoice {\hbox{\tenjap R}} {\hbox{\tenjap R}}
{\hbox{\sevenjap R}} {\hbox{\sevenjap R}}}}
\def\bbS{{\mathchoice{\hbox{\tenjap S}}{\hbox{\tenjap S}}
{\hbox{\sevenjap S}}{\hbox{\sevenjap S}}}}
\def\bbZ{{\mathchoice{\hbox{\tenjap Z}}{\hbox{\tenjap Z}}
{\hbox{\sevenjap Z}}{\hbox{\sevenjap Z}}}}

\def\w{{\mathchoice{\,{\scriptstyle\wedge}\,}
{{\scriptstyle\wedge}}
{{\scriptscriptstyle\wedge}}{{\scriptscriptstyle\wedge}}}}

\def\xb{{\bf x}}\def\yb{{\bf y}}\def\wb{{\bf w}}\def\zb{{\bf z}}
\def\ub{{\bf u}}\def\vb{{\bf v}}

\def\zerob{{\bf 0}}\def\oneb{{\bf 1}}
\def\Cl{{{\rm C}\ell}}\def\Im{{\rm Im}}
\def\bT{{\bf T}}\def\bR{{\bf R}}
\def\la{\langle}\def\ra{\rangle}
\def\euk{\hbox{\eu k}}\def\eur{\hbox{\eu r}}
\def\eug{\hbox{\eu g}}\def\euspin{\hbox{\eu spin}}\def\euso{\hbox{\eu so}}
\def\eusu{\hbox{\eu su}}
\def\eusl{\hbox{\eu sl}}
\def\GL{{\rm GL}}\def\SL{{\rm SL}}\def\End{\mathop{\rm End}}
\def\Spin{{\rm Spin}}\def\SO{{\rm SO}}\def\SU{{\rm SU}}
\def\Sp{{\rm Sp}}


\centerline{\bigbf Notes on Spinors in Low Dimension}
\bigskip
\centerline{Robert L. Bryant --- \tt{bryant@math.duke.edu} }
\bigskip
{\bf 0. Introduction.}  The purpose of these old notes (from 1998) is to determine the orbit structure of the groups~$\Spin(p,q)$ acting on their 
spinor spaces for certain values of~$n=p{+}q$, in particular, the values
$$
(p,q) = (8,0),\ (9,0),\ (10,0),\ \hbox{and}\ (10,1).
$$
though it will turn out in the end that there are a few interesting
things to say about the cases $(p,q) = (10,2)$ and $(9,1)$, as well.
\bigskip
{\bf 1. The Octonions.}  Let~$\bbO$ denote the ring of octonions. 
Elements of~$\bbO$ will be denoted by bold letters, such as $\xb$, $\yb$, 
etc.  Thus, $\bbO$ is the unique $\bbR$-algebra of dimension~$8$
with unit~$\oneb\in\bbO$ endowed with a positive definite inner 
product~$\la,\ra$ satisfying~$\la \xb\yb,\xb\yb\ra
=\la\xb,\xb\ra\,\la\yb,\yb\ra$ for all~$\xb,\yb\in\bbO$.   
As usual, the norm of
an element~$\xb\in\bbO$ is denoted $|\xb|$ and defined as the square root
of~$\la\xb,\xb\ra$.  Left and right multiplication
by~$\xb\in\bbO$ define maps~$L_\xb\,,R_\xb:\bbO\to\bbO$ that are isometries
when~$|\xb|=1$.

The conjugate of~$\xb\in\bbO$, 
denoted~$\overline\xb$, is defined to be~$\overline\xb = 
2\la\xb,\oneb\ra\,\oneb-\xb$.  When a symbol is needed, the map of 
conjugation will be denoted $C:\bbO\to\bbO$.  The identity
$\xb\,\overline\xb = |\xb|^2$ holds, as well as the conjugation 
identity~$\overline{\xb\yb}={\overline\yb}\,{\overline\xb}$.
In particular, this implies the useful identities 
$C\, L_\xb\, C=R_{\overline\xb}$
and $C\, R_\xb\, C=L_{\overline\xb}$.
\smallskip
The algebra~$\bbO$ is not commutative or associative.  However, any 
subalgebra of~$\bbO$ that is generated by two elements is associative.  
It follows that $\xb\,\bigl({\overline\xb} \yb\bigr)=|\xb|^2\,\yb$ and
that~$(\xb\yb)\xb=\xb(\yb\xb)$ for 
all~$\xb,\yb\in\bbO$. Thus,~$R_\xb\, L_\xb
=L_\xb\, R_\xb$ (though, of course, 
$R_\xb\, L_\yb\not=L_\yb\, R_\xb$ in general). In
particular, the expression $\xb\yb\xb$ is unambiguously defined.
In addition, there are the {\it Moufang Identities}
$$
\eqalign{
(\xb\yb\xb)\zb &= \xb\bigl(\yb(\xb\zb)\bigr),\cr
\zb(\xb\yb\xb) &= \bigl((\zb\xb)\yb\bigr)\xb,\cr
\xb(\yb\zb)\xb &= (\xb\yb)(\zb\xb),\cr
}
$$
which will be useful below. (See, for example,
{\it Spinors and Calibrations}, by F.~Reese Harvey, for proofs.)
\bigskip
{\bf 2.  $\Spin(8)$.}
For~$\xb\in \bbO$, define the linear 
map~$m_\xb:\bbO\oplus\bbO\to\bbO\oplus\bbO$ by the formula
$$
m_\xb = \left[\matrix{0&C \, R_\xb\cr -C \, L_\xb&0}\right]\,.
$$
By the above identities, it follows that $(m_\xb)^2 = -|\xb|^2$ and
hence this map induces a representation on the vector 
space~$\bbO\oplus\bbO$ of the Clifford algebra generated 
by~$\bbO$ with its standard quadratic form.  
This Clifford algebra is known to be isomorphic
to~$M_{16}(\bbR)$, the algebra of 16-by-16 matrices with real entries, so
this representation must be faithful.  By dimension count, this establishes
the isomorphism~$\Cl\bigl(\bbO,\la,\ra\bigr)
=\End_\bbR\bigl(\bbO\oplus\bbO\bigr)$.
\smallskip
The group $\Spin(8)\subset\GL_\bbR(\bbO\oplus\bbO)$ is defined as the
subgroup generated by products of the form~$m_\xb\,m_\yb$
where~$\xb,\yb\in\bbO$ satisfy~$|\xb|=|\yb|=1$.  
Such endomorphisms preserve the splitting of $\bbO\oplus\bbO$ into
the two given summands since
$$
m_\xb\,m_\yb 
= \left[\matrix{ -L_{\overline\xb} \, L_\yb&0\cr 
                 0&-R_{\overline\xb} \, R_\yb}\right]\,.
$$
In fact, setting $\xb=-\oneb$ in this formula shows that endomorphisms
of the form
$$
\left[\matrix{ L_\ub&0\cr 0&R_\ub}\right],\qquad\hbox{with $|\ub|=1$}
$$
lie in~$\Spin(8)$.  In fact, they generate~$\Spin(8)$, since 
$m_\xb\,m_\yb$ is clearly a product of two of these when 
$|\xb|=|\yb|=1$.
\smallskip
Fixing an identification~$\bbO\simeq\bbR^8$, this defines an embedding
$\Spin(8)\subset\SO(8)\times\SO(8)$, and the projections onto either
of the factors is a group homomorphism.  Since neither of these
projections is trivial, since the Lie algebra~$\euso(8)$ is simple, and
since~$\SO(8)$ is connected, it follows that each of these projections
is a surjective homomorphism.  Since~$\Spin(8)$ is simply connected
and since the fundamental group of~$\SO(8)$ is~$\bbZ_2$, it follows that  
that each of these homomorphisms is a non-trivial double cover of~$\SO(8)$.
Moreover, it follows that the subsets~$\{\  L_\ub\ \vrule\ |\ub|=1\ \}$
and $\{\  R_\ub\ \vrule\ |\ub|=1\ \}$ of~$\SO(8)$ each suffice to
generate~$\SO(8)$.
\smallskip
Let~$H\subset\bigl(\SO(8)\bigr)^3$ be the set of 
triples~$(g_1,g_2,g_3)\in\bigl(\SO(8)\bigr)^3$ for which
$$
g_2(\xb\yb) = g_1(\xb)\,g_3(\yb)
$$
for all~$\xb,\yb\in\bbO$.  The set~$H$ is closed and is evidently 
closed under multiplication and inverse. Hence it is a compact Lie group.  

By the third Moufang identity, $H$ contains the subset
$$
\Sigma = 
\left\{\ (L_\ub,\,L_\ub{\,}R_\ub,\,R_\ub)\ \vrule\ |\ub|=1
\right\}.  
$$
Let~$K\subset H$ be the subgroup generated by~$\Sigma$, and for $i=1,2,3$,
let~$\rho_i:H\to\SO(8)$ be the homomorphism that is 
projection onto the $i$-th factor.  Since~$\rho_1(K)$ contains
$\{\  L_\ub\ \vrule\ |\ub|=1\ \}$, it follows that $\rho_1(K)=\SO(8)$,
so {\it a fortiori} $\rho_1(H)=\SO(8)$.  Similarly, $\rho_3(H)=\SO(8)$.  

The kernel of~$\rho_1$ consists of elements~$(I_8,g_2,g_3)$
that satisfy~$g_2(\xb\yb) = \xb\,g_3(\yb)$ for all~$\xb,\yb\in\bbO$.  
Setting~$\xb=\oneb$ in this equation yields~$g_2=g_3$,
so that~$g_2(\xb\yb) = \xb\,g_2(\yb)$.  Setting~$\yb=\oneb$
in this equation yields~$g_2(\xb) = \xb\,g_2(\oneb)$, i.e., $g_2 = R_\ub$
for $\ub=g_2(\oneb)$.  Thus, the elements in the kernel of~$\rho_1$ are
of the form~$(1,R_\ub,R_\ub)$ for some~$\ub$ with~$|\ub|=1$.  
However, any such $\ub$ would, by definition, 
satisfy~$(\xb\yb)\ub=\xb(\yb\ub)$ for all~$\xb,\yb\in\bbO$, 
which is impossible unless~$\ub=\pm\oneb$.  Thus, the kernel of~$\rho_1$
is~$\bigl\{(I_8,\pm I_8,\pm I_8)\bigr\}\simeq\bbZ_2$, so that $\rho_1$ is a 
2-to-1 homomorphism of~$H$ onto~$\SO(8)$.  Similarly, $\rho_3$ is a 2-to-1 
homomorphism of~$H$ onto~$\SO(8)$, with 
kernel~$\bigl\{(\pm I_8,\pm I_8,I_8)\bigr\}$.  Thus, $H$ is either 
connected and isomorphic to~$\Spin(8)$ or else disconnected,
with two components.

Now~$K$ is a connected subgroup of~$H$ and the kernel of~$\rho_1$ 
intersected with~$K$ is either trivial or~$\bbZ_2$.  Moreover, the product 
homomorphism~$\rho_1{\times}\rho_3:K\to\SO(8){\times}\SO(8)$ maps the 
generator~$\Sigma\subset K$ into generators 
of~$\Spin(8)\subset\SO(8){\times}\SO(8)$.  It follows that
$\rho_1{\times}\rho_3(K)=\Spin(8)$ and hence that~$\rho_1$ and $\rho_3$ 
must be non-trivial double covers of~$\SO(8)$ when restricted to~$K$.  
In particular, it follows that~$K$ must be all of~$H$ and, moreover, that the 
homomorphism~$\rho_1{\times}\rho_3:H\to\Spin(8)$ must be an isomorphism.
It also follows that the homomorphism $\rho_2:H\to\SO(8)$ 
must be a double cover of~$\SO(8)$ as well.  

Henceforth, $H$ will be identified with~$\Spin(8)$ via the 
isomorphism~$\rho_1{\times}\rho_3$.  Note that the center of $H$
consists of the 
elements~$(\varepsilon_1\,I_8,\varepsilon_2\,I_8,\varepsilon_3\,I_8)$
where~${\varepsilon_i}^2=\varepsilon_1\varepsilon_2\varepsilon_3=1$
and is isomorphic to~$\bbZ_2\times\bbZ_2$.
\smallskip
{\it Triality.}  For~$(g_1,g_2,g_3)\in H$, the identity
$g_2(\xb\yb) = g_1(\xb)\,g_3(\yb)$ can be conjugated, giving
$$
Cg_2C(\xb\yb) 
= \overline{g_2(\overline{\yb}\,\overline{\xb})}
= \overline{g_1(\overline{\yb})\,g_3(\overline{\xb})}
= \overline{g_3(\overline{\xb})}\, \overline{g_1(\overline{\yb})}.
$$
This implies that $\bigl(Cg_3C,Cg_2C,Cg_1C\bigr)$ also lies in~$H$.
Also, replacing $\xb$ by $\zb\overline{\yb}$ in the original formula
and multiplying on the right by $\overline{g_3(\yb)}$ shows that
$$
g_2(\zb) \overline{g_3(\yb)} = g_1(\zb\overline{\yb}),
$$
implying that~$\bigl(g_2,g_1,Cg_3C\bigr)$ lies in~$H$ as well.  In fact, the
two maps~$\alpha,\beta:H\to H$ defined by
$$
\alpha(g_1,g_2,g_3) = \bigl(Cg_3C,Cg_2C,Cg_1C\bigr),
\quad\hbox{and}\qquad
\beta(g_1,g_2,g_3) = \bigl(g_2,g_1,Cg_3C\bigr) 
$$
are outer automorphisms (since they act nontrivially on the center of~$H$) 
and generate a group of automorphisms isomorphic to~$S_3$, the symmetric group 
on three letters.  The automorphism $\tau=\alpha\beta$ is known as the 
triality automorphism.
\smallskip
{\it Notation.} 
To emphasize the group action, denote~$\bbO\simeq\bbR^8$ by~$V_i$
when regarding it as a representation space of~$\Spin(8)$ via
the representation~$\rho_i$.  Thus, octonion multiplication induces
a $\Spin(8)$-equivariant projection
$$
V_1\otimes V_3 \longrightarrow V_2\,.
$$
In the standard notation, it is traditional to identify $V_1$ with
$\bbS_-$ and~$V_3$ with $\bbS_+$ and to refer to $V_2$ as the 
`vector representation'.  Let~$\rho_i':\euspin(8)\to\euso(8)$ denote
the corresponding Lie algebra homomorphisms, which are, in fact, 
isomorphisms. For simplicity of notation, for any~$a\in\euspin(8)$,
the element~$\rho_i'(a)\in\euso(8)$ will be denoted by~$a_i$ when
no confusion can arise.
\smallskip
{\it Orbit structure.}  Let~$\SO(\Im\bbO)\simeq\SO(7)$ denote the 
subgroup of~$\SO(\bbO)\simeq\SO(8)$ that leaves~$\oneb\in\bbO$ fixed,
and let~$K_i\subset H$ be the preimage of~$\SO(\Im\bbO)$ under the
homomorphism~$\rho_i:H\to\SO(\bbO)$.  Then~$K_i$ is a non-trivial
double cover of~$\SO(\Im\bbO)$ and hence is isomorphic to~$\Spin(7)$.
Note, in particular that~$K_1$ contains~$(I_8,-I_8,-I_8)$ and hence
$\rho_3(K_1)\subset\SO(8)$ contains~$-I_8$.  This implies that
$\rho_3:K_1\to\SO(V_3)$ is a faithful representation of~$\Spin(7)$ and
hence~$K_1$ acts transitively on the unit sphere in~$V_3$.  
\smallskip
In particular, it follows that~$\Spin(8)\subset\SO(V_1)\times\SO(V_3)$
acts transitively on the product of the unit spheres in~$V_1$ and~$V_3$.
Consequently, it follows that the quadratic polynomials
$$
q_1(\xb,\yb) = |\xb|^2\qquad\rm{and}\qquad q_2(\xb,\yb) = |\yb|^2
$$
generate the ring of $\Spin(8)$-invariant polynomials on~$\bbO\oplus\bbO$
and that every point of this space lies on the $\Spin(8)$-orbit of a
unique element of the form~$(a\,\oneb,b\,\oneb)$ for some pair of
real numbers~$a,b\ge0$.  For~$ab\not=0$, the stabilizer of such an element
is the 14-dimensional simple group~$G_2$, and this group acts transitively
on the unit sphere in~$\Im\bbO$.
\bigskip
{\bf 3. $\Spin(9)$.}  For~$(r,\xb)\in\bbR\oplus\bbO$, define a
$\bbC$-linear map~$m_{(r,\xb)}:\bbC\otimes\bbO^2\to\bbC\otimes\bbO^2$
by the formula
$$
m_{(r,\xb)}
=i \left[\matrix{ r\,I_8& C\,R_\xb\cr C\,L_\xb& -r\,I_8}\right]\,.
$$
Since~$(m_{(r,\xb)})^2$ is $-\left(r^2{+}|\xb|^2\right)$ times the identity 
map, this defines a $\bbC$-linear representation 
on~$\bbC\otimes\bbO^2$ of the Clifford algebra generated by~$\bbR\oplus\bbO$ 
endowed with its direct sum inner product.  Since this Clifford algebra is 
known to be isomorphic to~$M_{16}(\bbC)$, it follows, for dimension reasons, 
that this representation is one-to-one and onto, establishing the 
isomorphism~$\Cl(\bbR\oplus\bbO,\la,\ra) = \End_\bbC(\bbC\otimes\bbO^2)$.

As usual, $\Spin(9)$ is the subgroup generated by the products of the
form~$m_{(r,\xb)}m_{(s,\yb)}$ where $r^2+|\xb|^2=s^2+|\yb|^2=1$. Note
that these products have real coefficients, and so actually lie in
$\GL_\bbR(\bbO^2)\simeq\GL(16,\bbR)$.  In fact, these products are
themselves seen to be products of the products of the special form
$$
p_{(r,\xb)}  
= m_{(-1,\zerob)}m_{(r,\xb)} 
= \left[\matrix{ r\,I_8& C\,R_\xb\cr -C\,L_\xb&  r\,I_8}\right]\,,
\qquad\hbox{where $r^2+|\xb|^2=1$},
$$
so these latter matrices suffice to generate~$\Spin(9)$. By the results
of the previous section, products of an even number of the~$p_{(0,\ub)}$
with~$|\ub|=1$
generate~$\Spin(8)\subset\Spin(9)$.  
\smallskip
Since the linear transformations of the form~$p_{(r,\xb)}$ preserve the
quadratic form
$$
q(x,y) = |\xb|^2+|\yb|^2,
$$
it follows that~$\Spin(9)$ is a subgroup of~$\SO(\bbO^2)=\SO(16)$.
\medskip
{\it The Lie algebra.}  Since~$\Spin(9)$ contains $\Spin(8)$, the
containment~$\euspin(8)\subset\euspin(9)$ yields the containment
$$
\left\{\ \pmatrix{a_1&0\cr0&a_3}\ \vrule\ a\in\euspin(8)\ \right\}
\subset\euspin(9)\,.
$$
Moreover, since~$\Spin(9)$ contains the 8-sphere consisting of the
$p_{(r,\xb)}$ with $r^2+|\xb|^2=1$, its Lie algebra must contain
the tangent space to this 8-sphere at $(r,\xb)=(1,\zerob)$, i.e.,
$$
\left\{\ \pmatrix{0&C\,R_\xb\cr-C\,L_\xb&0}\ \vrule\ \xb\in\bbO\ \right\}
\subset\euspin(9)\,.
$$
By dimension count, this implies the equality
$$
\euspin(9) = \left\{\pmatrix{a_1&C\,R_\xb\cr-C\,L_\xb&a_3}\ 
  \vrule\ \xb\in\bbO,\ a\in\euspin(8)\ \right\}\,.
$$

Let~$\rho:\Spin(9)\to\SO(\bbR\oplus\bbO)\simeq\SO(9)$ be the
homomorphism for which the induced map on Lie algebras is
$$
\rho'\left(\pmatrix{a_1&C\,R_\zb\cr-C\,L_\zb&a_3}\right)
= \pmatrix{0&2\,\overline{\zb}^*\cr-2\,\overline{\zb}&a_2}.
$$
where~$\xb^*:\bbO\to\bbR$ is just~$\xb^*(\yb) = \la\xb,\yb\ra$.
(The triality constructions imply that $\rho'$ is, indeed, a Lie algebra 
homomorphism.  Note that, when restricted to~$\Spin(8)$, this becomes
the homomorphism~$\rho_2:\Spin(8)\to\SO(\bbO)=\SO(8)$.)  Then~$\rho$
is a double cover of~$\SO(9)$.  
\smallskip
Define the squaring map~$\sigma:\bbO^2\to\bbR\oplus\bbO$ by
$$
\sigma\left(\pmatrix{\xb\cr\yb\cr}\right) 
= \pmatrix{|\xb|^2-|\yb|^2\cr 2\,\xb\,\yb\cr}.
$$
A short calculation using the Moufang Identities shows that $\sigma$ is
$\rho$-equivariant, i.e., that~$\sigma\bigl(g\vb\bigr) 
=\rho(g)\bigl(\sigma(\vb)\bigr)$ for all~$\vb\in\bbO^2$ 
and all~$g\in\Spin(9)$.  This will be useful below.
\medskip 
{\it Orbit structure and stabilizer.}  Each point of~$\bbO^2$ lies
on the $\Spin(8)$-orbit of an element~$(a\,\oneb,b\,\oneb)$ 
for some real numbers~$a,b\ge0$. Thus, the orbits of~$\Spin(9)$ on 
the unit sphere in~$\bbO^2$ are unions of the $\Spin(8)$-orbits of the 
elements~$(\cos\theta\,\oneb,\sin\theta\,\oneb)$.  
Now, calculation yields
$$
p_{(\cos\phi,\sin\phi\,\oneb)}\pmatrix{\cos\theta\,\oneb\cr\sin\theta\,\oneb}
= \pmatrix{\cos(\theta{-}\phi)\,\oneb\cr\sin(\theta{-}\phi)\,\oneb}.
$$
Since all of the elements~$(\cos\theta\,\oneb,\sin\theta\,\oneb)$ lie on
a single $\Spin(9)$-orbit, it follows that~$\Spin(9)$ acts transitively
on the unit sphere in~$\bbO^2$ and, consequently, that the quadratic 
form~$q$ generates the ring of $\Spin(9)$-invariant polynomials 
on~$\bbO^2$.

Since the orbit of~$(\oneb,\zerob)\in\bbO^2$ is the 15-sphere and
since $\Spin(9)$ is connected and simply connected, it follows that
the $\Spin(9)$-stabilizer of this element must be connected, simply connected,
and of dimension~$21$.  Since~$K_1\subset\Spin(8)\subset\Spin(9)$
lies in this stabilizer and has dimension~$21$, it follows that~$K_1$
must be equal to this stabilizer.
\smallskip
 For use in the next two sections, it will be useful to understand the
orbits of~$\Spin(9)$ acting on~$\bbO^2\oplus\bbO^2$ and to understand
the ring of $\Spin(9)$-invariant polynomials on this vector space of
real dimension~$32$.  The first observation is that the generic orbit
has codimension~$4$.  This can be seen as follows:  Since~$\Spin(9)$
acts transitively on the unit sphere in~$\bbO^2$, every element lies
on the ~$\Spin(9)$ orbit of an element of the form
$$
\left(\ \pmatrix{a\,\oneb\cr\zerob\cr},
      \ \pmatrix{\xb\cr\yb\cr}\ \right),
$$
where~$a\ge0$.  Assuming~$a>0$, the stabilizer in~$\Spin(9)$ of this
first component is~$K_1\simeq\Spin(7)$ and this acts transitively on the
unit sphere in the second $\bbO$-summand of~$\bbO^2$, so that an
element of the above form lies on the orbit of an element of the form
$$
\left(\ \pmatrix{a\,\oneb\cr\zerob\cr},
      \ \pmatrix{\xb\cr b\,\oneb\cr}\ \right),
$$
where~$b\ge0$.  Assuming~$b>0$, the stabilizer in~$K_1$ of $\oneb$ in
this second $\bbO$-summand is~$G_2$, which acts transitively on the
unit sphere in~$\Im\bbO$ in the first $\bbO$-summand.  This implies that
an element of the above form lies on the orbit of an element of the form
$$
\zb = \left(\ \pmatrix{a\,\oneb\cr\zerob\cr},
      \ \pmatrix{c\,\oneb + d\,\ub \cr b\,\oneb\cr}\ \right),
$$
for some~$c,d\ge0$ and~$\ub\in\Im\bbO$ some fixed unit imaginary octonion.
Thus, the generic $\Spin(9)$-orbit has codimension at most~$4$.  It is
still possible that two elements of the above form with distinct
values of~$a,b,c,d>0$ might lie on the same~$\Spin(9)$-orbit, but this
will be ruled out directly.

To see that these latter elements lie on distinct $\Spin(9)$-orbits, 
it will be sufficient to construct $\Spin(9)$-invariant polynomials 
on~$\bbO^2\oplus\bbO^2$ that separate these elements.  To do so,
write the typical element of~$\bbO^2\oplus\bbO^2$ in the form
$$
(\vb_1,\vb_2) = \left(\ \pmatrix{\xb_1\cr\yb_1\cr},
      \ \pmatrix{\xb_2\cr\yb_2\cr}\ \right),
$$
and first consider the three quadratic polynomials
$$
\eqalign{
q_{2,0} &= |\xb_1|^2+|\yb_1|^2\cr
q_{1,1} &= \xb_1\cdot\xb_2+\yb_1\cdot\yb_2\cr
q_{0,2} &= |\xb_2|^2+|\yb_2|^2\cr
}
$$
These polynomials are manifestly $\Spin(9)$-invariant and satisfy
$$
q_{2,0}(\zb) = a^2,\qquad q_{1,1}(\zb) = ac,
\qquad q_{0,2}(\zb) = b^2+ c^2 + d^2.
$$
Evidently, these polynomials span the vector
space of $\Spin(9)$-invariant quadratic polynomials on~$\bbO^2\oplus\bbO^2$.

Since $\Spin(9)$ contains $-1$ times the identity, there are no 
$\Spin(9)$-invariant cubic polynomials.  A representation-theoretic argument
shows that the $\Spin(9)$-invariant quartic polynomials 
on~$\bbO^2\oplus\bbO^2$ form a vector space of dimension~$7$.  Six of these 
are accounted for by quadratic polynomials in~$q_{2,0}$, $q_{1,1}$, 
and $q_{0,2}$, while a seventh can be constructed as follows.  Define 
$$
q_{2,2} = \sigma(\vb_1)\cdot\sigma(\vb_2)
        = \left(|\xb_1|^2-|\yb_1|^2\right)\left(|\xb_2|^2-|\yb_2|^2\right)
          + 4\,(\xb_1\yb_1)\cdot(\xb_2\yb_2).
$$
Using the $\Spin(9)$-equivariance of the squaring map~$\sigma$, it follows
that $q_{2,2}$ is indeed invariant under~$\Spin(9)$.
Note that
$$
q_{2,2}(\zb) = a^2(c^2+d^2-b^2),
$$
so that knowledge 
of~$\bigl(q_{2,0}(\zb),q_{1,1}(\zb),q_{0,2}(\zb),q_{2,2}(\zb) \bigr)$ 
suffices to recover~$a,b,c,d>0$ when these numbers are all non-zero. It 
now follows that the simultaneous level sets of these four polynomials are 
exactly the $\Spin(9)$-orbits on~$\bbO^2\oplus\bbO^2$.  (It seems likely 
that these polynomials generate the ring of $\Spin(9)$-invariant polynomials
on~$\bbO^2\oplus\bbO^2$, but such a result will not be needed, so this
problem will not be discussed further.)
\bigskip
{\bf 4. $\Spin(10)$.}  Rather than construct the Clifford representation
for an inner product on a vector space of dimension~$10$, it is convenient
to use the fact that~$\Spin(10)$ already appears as a subgroup 
of~$\Cl(\bbR\oplus\bbO,\la,\ra)=\End_\bbC(\bbC\otimes\bbO^2)$.  In fact, by 
the discussion in the last section, $\Spin(10)$ is the connected subgroup of 
this latter algebra whose Lie algebra is
$$
\euspin(10) 
= \left\{\pmatrix{a_1+ir\,I_8&C\,R_\xb+i\,C\,R_\yb\cr
                    -C\,L_\xb+i\,C\,L_\yb&a_3-ir\,I_8}\ 
  \vrule\ r\in\bbR,\ \xb,\yb\in\bbO,\ a\in\euspin(8)\ \right\}\,.
$$
Note that~$\euspin(10)$ appears as a subspace of~$\eusu(16)$, so that
$\Spin(10)$ acts on~$\bbC^{16}=\bbC\otimes\bbO^2$ preserving the complex
structure and the quadratic form
$$
q = q_{2,0}+q_{0,2} = |\xb_1|^2+|\yb_1|^2 + |\xb_2|^2+|\yb_2|^2\,,
$$
where, now, the typical element of~$\bbC\otimes\bbO^2$ will be written
as
$$
\zb = \pmatrix{\xb_1+i\,\xb_2\cr \yb_1+i\,\yb_2\cr}.
$$
\smallskip
Note that, because there are no connected Lie groups that lie properly
between~$\Spin(9)$ and $\Spin(10)$, it follows that~$\Spin(10)$ is
generated by~$\Spin(9)$ and the circle subgroup
$$
\bT 
=\left\{\pmatrix{e^{ir}\,I_8&0\cr0&e^{-ir}\,I_8}
  \ \vrule\ r\in\bbR/2\pi\bbZ\ \right\},
$$
which lies in~$\Spin(10)$, but does not lie in~$\Spin(9)$.  In particular,
a polynomial on~$\bbC\otimes\bbO^2$ is~$\Spin(10)$-invariant if and only
if it is both~$\Spin(9)$-invariant and~$\bT$-invariant.
\medskip
{\it Invariant polynomials.}
Among the quadratic polynomials that are $\Spin(9)$-invariant, only the 
multiples of~$q = q_{2,0}+q_{0,2}$ are also~$\bT$-invariant.  Thus, $q$ spans
the space of $\Spin(10)$-invariant quadratic forms on~$\bbC\otimes\bbO^2$.
In particular, this implies that the action of~$\Spin(10)$ 
on~$\bbC\otimes\bbO^2$ is irreducible (even as a real vector space).  

Among the quartic polynomials that are~$\Spin(9)$-invariant, a short 
calculation shows that only linear combinations of~$q^2$ and 
$$
\eqalign{
p &= {\ts{1\over2}}\left(q_{2,2} + q_{2,0}\,q_{0,2}-2\,{q_{1,1}}^2\right)\cr
  &= |\xb_1|^2|\xb_2|^2+|\yb_1|^2|\yb_2|^2 
     - \left(\xb_1\cdot\xb_2+\yb_1\cdot\yb_2\right)^2
     + 2\,(\xb_1\yb_1)\cdot(\xb_2\yb_2)\cr
  &= |\xb_1\w\xb_2|^2 + |\yb_1\w\yb_2|^2     
     - 2\,(\xb_1\cdot\xb_2)\,(\yb_1\cdot\yb_2)
     + 2\,(\xb_1\yb_1)\cdot(\xb_2\yb_2) .
}
$$
are invariant under the action of~$\bT$.  Thus, it follows that $q^2$~and~$p$ 
span the space of $\Spin(10)$-invariant quartics.  (Note the interesting
feature that, in the latter expression for~$p$, only the final term makes
use of octonion multiplication operations.)
\medskip
{\it Orbits and stabilizers.}  Let~$M\subset\bbC\otimes\bbO^2$ be the
$\Spin(10)$-orbit of~$\zb_0=(\oneb+i\,\zerob,\zerob + i\,\zerob)$.  The
tangent space to~$M$ at~$\zb_0$ is the set of vectors of the form
$$
\pmatrix{a_1+ir\,I_8&C\,R_\xb+i\,C\,R_\yb\cr
                    -C\,L_\xb+i\,C\,L_\yb&a_3-ir\,I_8}
\pmatrix{\oneb+i\,\zerob\cr\zerob + i\,\zerob}
= \pmatrix{a_1\oneb+i\,r\oneb\cr -\overline{\xb} + i\,\overline{\yb}}\,.
$$ 
and the Lie algebra of the $\Spin(10)$-stabilizer of~$\zb_0$ is defined
by the equations~$a_1\oneb = r = \xb = \yb = 0$.  Thus, the
identity component of the stabilizer of~$\zb_0$ is~$K_1\simeq\Spin(7)$
and the full stabilizer must lie in the normalizer of~$K_1$ in~$\Spin(10)$.
Evidently, the normalizer of~$K_1$ in~$\Spin(10)$ is~$K_1\cdot\bT$.
Since only the identity in the subgroup~$\bT$ stabilizes~$\zb_0$, 
the full stabilizer of~$\zb_0$ is~$K_1$.   Thus, $M$ is
diffeomorphic to~$\Spin(10)/\Spin(7)$, which is a smooth manifold of 
dimension~$45-21=24$ that is 2-connected, i.e., $\pi_0(M)=\pi_1(M)=\pi_2(M)=0$.

The normal space to~$M$ at~$\zb_0$ is the orthogonal direct sum of the
line~$\bbR\zb_0$ (which is normal to the unit sphere in~$\bbC\otimes\bbO^2$)
and the subspace of dimension~$7$
$$
N_{\zb_0} = 
\left\{\pmatrix{0+i\,\xb\cr\zerob+i\,\zerob}\ \vrule\ \xb\in\Im\bbO\ \right\}.
$$
The stabilizer~$K_1$ acts as~$\SO(7)$ on this
subspace.  In particular, it acts transitively on the unit sphere 
in~$N_{\zb_0}$, and hence it acts transitively on the space of geodesics
in the unit 31-sphere that meet~$M$ orthogonally at~$\zb_0$.  Since~$M$
is itself a $\Spin(10)$-orbit, it follows that~$\Spin(10)$ must 
act transitively on the the normal tube of any radius about~$M$ in the
unit 31-sphere.  Since, for generic radii, these normal tubes are 
hypersurfaces, it follows that the generic $\Spin(10)$-orbit in the 
31-sphere must be a hypersurface of dimension~$30$.   
Using the fact that such a hypersurface is an $S^6$-bundle over~$M$,
the long exact sequence in homotopy implies that these hypersurface
orbits are also 2-connected, which implies that the $\Spin(10)$-stabilizer
of any point on such a hypersurface must be both connected and
simply connected.

Now, the full group~$\Spin(10)$ must act transitively on the
space of geodesics in the unit 31-sphere that meet~$M$ orthogonally at
any point while every point of the unit 31-sphere lies on some
geodesic that meets~$M$ orthogonally. Thus, fixing some $\ub\in\Im\bbO$ 
with~$|\ub|=1$, it follows that every element of the 31-sphere lies on 
the $\Spin(10)$-orbit of an element of the form
$$
\zb_\theta = \pmatrix{\cos\theta+i\,\sin\theta\,\ub\cr\zerob+i\,\zerob}.
$$
Note that $p(\zb_\theta) = \cos^2\theta\,\sin^2\theta 
= {\ts{1\over4}}\,\sin^2(2\theta)$, so it follows that for~$0\le \theta
\le \pi/4$, the elements~$\zb_\theta$ lie on distinct orbits, and that
$0\le p\le {\ts{1\over4}}$, with the endpoints of this interval being
the only critical values of~$p$.  While~$M = p^{-1}(0)$ is one critical
orbit, the other critical orbit is~$M^* = p^{-1}({1\over4})$ and consists
of the points of the 31-sphere that are at geodesic 
distance~$\pi/4$ from~$M$.  It follows from this that~$M^*$ is also
connected and is a single orbit of~$\Spin(10)$.  In particular, the
simultaneous level sets of~$q$ and $p$ are exactly the~$\Spin(10)$-orbits
in~$\bbC\otimes\bbO^2$.

 For~$0<\theta<\pi/4$, the nearest point on~$M$ to~$\zb_\theta$ is~$\zb_0$,
so the $\Spin(10)$-stabilizer of~$\zb_\theta$ is a subgroup of~$K_1$
that has already been seen to be both connected and simply connected. 
Also, the orbit of~$\zb_\theta$ is a 6-sphere bundle over~$M$.  By
dimension count, this stabilizer must be of dimension 15 and must contain 
the stabilizer in~$K_1$ of~$\oneb$ and $\ub$, which is~$\Spin(6)$.  
Thus, the stabilizer of such a~$\zb_\theta$ is exactly $\Spin(6)\simeq\SU(4)$. 
In particular, the stabilizer of any point of the 31-sphere not 
on~$M$ or~$M^*$ must be a conjugate of~$\SU(4)$.

Now, the tangent space to~$M^*$ at~$\zb_{\pi/4}$ is the set of vectors
of the form
$$
\pmatrix{a_1+ir\,I_8&C\,R_\xb+i\,C\,R_\yb\cr
                    -C\,L_\xb+i\,C\,L_\yb&a_3-ir\,I_8}
\pmatrix{\oneb+i\,\ub\cr\zerob + i\,\zerob}
= 
  \pmatrix{(a_1\oneb-r\,\ub)+i\,(a_1\ub+r\oneb)\cr\noalign{\vskip2pt}
           -\overline{(\xb+\yb\ub)}+i\,\overline{(\yb-\xb\ub)}}\,.
$$ 
Thus, the Lie algebra of the stabilizer~$G$ of~$\zb_{\pi/4}$ is
defined by the relations~$a_1\oneb-r\,\ub=a_1\ub+r\oneb=\yb-\xb\ub=\zerob$.
(Remember that~$\ub^2 = -\oneb$.)  It follows that~$a_1\in\euso(8)$ must
belong to the stabilizer of the 2-plane spanned by~$\{\oneb,\ub\}$, so
that~$a_1$ lies in~$\euso(2)\oplus\euso(6)$.  Conversely, if~$a_1$ lies
in this subspace, then there exists a unique~$r\in\bbR$ so that 
$a_1\oneb-r\,\ub=a_1\ub+r\oneb = 0$.  From the matrix representation, it is
clear that the maximal torus in $\euso(2)\oplus\euso(6)$ (which has rank~4)
is a maximal torus in the full stabilizer algebra, which has dimension~$24$.
The root pattern is evident from the matrix representation, implying
that the stabilizer algebra is isomorphic to~$\eusu(5)$.  

Now~$M^*$ has dimension~$21$ and is the base of a fibration whose total space 
is one of the hypersurface orbits and whose fiber is a 9-sphere.  
The 2-connectivity of the hypersurface orbits implies, by the long
exact sequence in homotopy, that~$M^*$ is also 2-connected,  which implies
that $M^* = \Spin(10)/G$ where $G$ is both connected and simply connected.
Since its Lie algebra is~$\eusu(5)$, it follows that~$G$ is isomorphic 
to~$\SU(5)$.
\bigskip
{\bf 5. $\Spin(10,1)$.}  To construct the spinor representation 
of~$\Spin(10,1)$, it will be easiest to construct the Lie algebra 
representation by extending the Lie algebra representation of~$\Spin(10)$
that was constructed in~\S3.  It is convenient to identify~$\bbC\otimes\bbO^2$
with~$\bbO^4$ explicitly via the identification
$$
\zb = \pmatrix{\xb_1+i\,\xb_2\cr \yb_1+i\,\yb_2\cr}
    = \pmatrix{\xb_1\cr\yb_1\cr\xb_2\cr\yb_2\cr}.
$$
Via this identification, $\euspin(10)$ becomes the subspace
$$
\euspin(10) 
= \left\{\pmatrix{
   a_1    & C\,R_\xb &  -r\,I_8  &  -C\,R_\yb\cr
-C\,L_\xb & a_3      & -C\,L_\yb &  r\,I_8   \cr 
r\,I_8    & C\,R_\yb &    a_1    &   C\,R_\xb\cr
 C\,L_\yb &  -r\,I_8 & -C\,L_\xb &     a_3   \cr}\ 
  \vrule\ r\in\bbR,\ \xb,\yb\in\bbO,\ a\in\euspin(8)\ \right\}\,.
$$
Consider the one-parameter subgroup~$\bR\subset\SL_\bbR(\bbO^4)$ defined by
$$
\bR 
=\left\{\pmatrix{t\,I_{16}&0\cr0&t^{-1}\,I_{16}}
  \ \vrule\ t\in\bbR^+\ \right\}.
$$
It has a Lie algebra~$\eur\subset\eusl(\bbO^4)$. Evidently, the
the subspace~$\left[\euspin(10),\eur\right]$ consists of matrices of
the form
$$
\pmatrix{
   0_8 & 0_8 &  r\,I_8  &  C\,R_\yb\cr
   0_8 & 0_8 & C\,L_\yb &  -r\,I_8   \cr 
 r\,I_8    & C\,R_\yb & 0_8 & 0_8 \cr
 C\,L_\yb &  -r\,I_8 & 0_8 &  0_8 \cr},
\qquad r\in\bbR,\ \yb\in\bbO\,.
$$

Let~$\eug = \euspin(10) \oplus\eur\oplus\left[\euspin(10),\eur\right]$. 
Explicitly,
$$
\eug 
= \left\{\pmatrix{
   a_1 +x\,I_8  & C\,R_\xb &  y\,I_8  &  C\,R_\yb\cr
-C\,L_\xb & a_3 +x\,I_8  & C\,L_\yb &  -y\,I_8   \cr 
z\,I_8    & C\,R_\zb &  a_1 -x\,I_8   &   C\,R_\xb\cr
 C\,L_\zb &  -z\,I_8 & -C\,L_\xb &     a_3 -x\,I_8   \cr}\ 
  \vrule\ 
\matrix{x,y,z\in\bbR,\cr\noalign{\vskip2pt} \xb,\yb,\zb\in\bbO,\cr
        \noalign{\vskip2pt}\ a\in\euspin(8)}\ \right\}\,.
$$
Compuation using the Moufang Identities shows that $\eug$ is closed
under Lie bracket and hence is a Lie algebra of (real) dimension~$55$
that contains~$\euspin(10)$.  The induced representation 
of~$\Spin(10)$ on~$\eug/\euspin(10)$ evidently restricts to~$\Spin(9)$
to preserve the splitting corresponding to the 
sum~$\eur\oplus\left[\euspin(10),\eur\right]\simeq\bbR\oplus\bbR^9$ and
acts as the standard (irreducible) representation on the ~$\bbR^9$ summand.
It follows that~$\Spin(10)$ must act via its standard (irreducible, ten
dimensional) representation on~$\eug/\euspin(10)$.  Since the trace of the 
square of a non-zero element in the 
subspace~$\eur\oplus\left[\euspin(10),\eur\right]$ is positive, $\eug$ is
semisimple of non-compact type. It follows that $\eug$ is isomorphic 
to~$\euso(10,1)$ and hence is the Lie algebra of a representation 
of~$\Spin(10,1)$.  This representation must be faithful since it is faithful 
on the maximal compact subgroup~$\Spin(10)$.

Thus, define~$\Spin(10,1)$ to be the (connected) subgroup 
of~$\SL_\bbR(\bbO^4)$ that is generated by~$\Spin(10)$ and the 
subgroup~$\bR$.  Its Lie algebra~$\eug$ will henceforth be written 
as~$\euspin(10,1)$.
\medskip
{\it Invariant Polynomials and Orbits.}  Consider the 
$\Spin(10)$-invariant polynomial
$$
p = |\xb_1|^2|\xb_2|^2+|\yb_1|^2|\yb_2|^2 
     - \left(\xb_1\cdot\xb_2+\yb_1\cdot\yb_2\right)^2
     + 2\,(\xb_1\yb_1)\cdot(\xb_2\yb_2)\,,
$$
Evidently, $p$ is invariant under~$\bR$ and is therefore invariant under
$\Spin(10,1)$.  In particular, it follows that the orbits of~$\Spin(10,1)$
on~$\bbO^4\simeq\bbR^{32}$ must lie in the level sets of~$p$.  

Also from the previous section, it is known that every element of~$\bbO^4$
lies on the $\Spin(10)$-orbit of exactly one of the elements
$$
\zb_{a,b} = \pmatrix{a\,\oneb\cr\zerob\cr b\,\ub\cr\zerob}
\qquad\hbox{ where $0\le b\le a$.}
$$
and where~$\ub\in\Im\bbO$ is a fixed unit imaginary octonion. 
However, all of the elements of the form
$$
\pmatrix{at\,\oneb\cr\zerob\cr (b/t)\,\ub\cr\zerob}
\qquad\hbox{(where $0< t$)}
$$
lie on the same~$\bR$-orbit and, hence, on the same
$\Spin(10,1)$-orbit.  Since $p(\zb_{a,b}) = a^2b^2$, it now follows that
each of the nonzero level sets of~$p$ is a single $\Spin(10,1)$-orbit while
the zero level set is the union of the origin and a single 
$\Spin(10,1)$-orbit, say, the orbit of~$\zb_{1,0}$.  
Moreover, it follows that $p$ generates the ring of $\Spin(10,1)$-invariant 
polynomials on~$\bbO^4$.
\medskip
{\it Stabilizers.}  Multiplication by positive scalars acts transitively
on the non-zero level sets of~$p$, so they are all diffeomorphic.  In
fact, each such level set is contractible to the $\Spin(10)$-invariant
locus where~$q$ reaches its minimum on this level set and this is a
manifold of dimension~$21$ that is diffeomorphic to~$M^*$.  In particular,
it follows that each of the non-zero level sets of~$p$ is 2-connected,
so that the stabilizer in~$\Spin(10,1)$ of a point on such a level set
must be connected and simply connected. 

If~$\zb\in\bbO^4$ has $p(\zb)\not=0$, then the $\Spin(10,1)$-orbit of~$\zb$ 
has dimension 31 and so its stabilizer in~$\Spin(10,1)$ must be of 
dimension~$55{-}31=24$. Moreover all of these stabilizers must be
conjugate in~$\Spin(10,1)$.  Since the $\Spin(10)$-stabilizer of the
point~$\zb_{1,1}$ is already known to be~$\SU(5)$, which has dimension~$24$,
it follows that this must be the~$\Spin(10,1)$-stabilizer as well.

The $\Spin(10,1)$-orbit consisting of nonzero vectors in the zero locus 
of~$p$ is just the deleted cone on~$M$, and so has dimension~$25$.  Since it 
is contractible to~$M$, it is 2-connected, so that the stabilizer 
in~$\Spin(10,1)$ of a point in this orbit must be connected and simply
connected and of dimension~$55{-}25=30$.  In fact, the Lie algebra of
this stabilizer is just
$$ 
\left\{\pmatrix{
a_1& 0 &  y\,I_8  &  C\,R_\yb\cr
 0 &a_3& C\,L_\yb &  -y\,I_8   \cr 
 0 & 0 &a_1& 0 \cr
 0 & 0 & 0 &a_3\cr}\ 
  \vrule\ 
\matrix{y\in\bbR,\cr\noalign{\vskip2pt} \yb\in\bbO,\cr
        \noalign{\vskip2pt}\ a\in\euk_1}\ \right\}
$$
where~$\euk_1$ is the Lie algebra of~$K_1\subset\Spin(8)$.  Thus, the
stabilizer is a semi-direct product of~$\Spin(7)$ with a copy of~$\bbR^9$.
\bigskip
Consider the squaring 
map~$\sigma:\bbO^4\to\bbR^{2+1}\oplus\bbO=\bbR^{10+1}$ that takes spinors 
for~$\Spin(10,1)$ to vectors.  This map~$\sigma$ is defined as follows:
$$
\sigma\left(\pmatrix{\xb_1\cr\yb_1\cr\xb_2\cr\yb_2\cr}\right) 
= \pmatrix{|\xb_1|^2+|\yb_1|^2\cr 
        2\,\bigl(\xb_1\cdot\xb_2-\yb_1\cdot\yb_2\bigr)\cr
        |\xb_2|^2+|\yb_2|^2\cr
      2\,\bigl(\xb_1\,\yb_2+\xb_2\,\yb_1\bigr)\cr}.
$$
Define the inner product on vectors in~$\bbR^{2+1}\oplus\bbO$ by the rule
$$
\pmatrix{a_1\cr a_2\cr a_3\cr \xb}\cdot
\pmatrix{b_1\cr b_2\cr b_3\cr \yb} 
= -2(a_1b_3+a_3b_1)+a_2b_2 + \xb\cdot\yb
$$
and let~$\SO(10,1)$ denote the subgroup of~$\SL(\bbR^{2+1}\oplus\bbO)$
that preserves this inner product.  This group still has two components
of course, but only the identity component~$\SO^\uparrow(10,1)$ will be of 
interest here. Let~$\rho:\Spin(10,1)\to\SO^\uparrow(10,1)$ be the homomorphism 
whose induced map on Lie algebras is given by the isomorphism
$$
\rho'\left(
\pmatrix{
   a_1 +x\,I_8  & C\,R_\xb &  y\,I_8  &  C\,R_\yb\cr
-C\,L_\xb & a_3 +x\,I_8  & C\,L_\yb &  -y\,I_8   \cr 
z\,I_8    & C\,R_\zb &  a_1 -x\,I_8   &   C\,R_\xb\cr
 C\,L_\zb &  -z\,I_8 & -C\,L_\xb &     a_3 -x\,I_8   \cr}
\right)
= 
\pmatrix{
   2x  & y &  0  &  \overline{\yb}^*\cr
2z & 0  & 2y &  2\,\overline{\xb}^*   \cr 
0   & z &   -2x   &   \overline{\zb}^*\cr
 2\,\overline{\zb} &  -2\,\overline{\xb} & 2\,\overline{\yb} & a_2\cr}.
$$
With these definitions, the squaring map~$\sigma$ is seen to have the
equivariance~$\sigma\bigl(g\,\zb\bigr) = \rho(g)\,\bigl(\sigma(\zb)\bigr)$ 
for all~$g$ in~$\Spin(10,1)$ and all~$\zb\in\bbO^4$.  

With these definitions, the polynomial~$p$ has the expression~$p(\zb) 
= -{\ts{1\over4}}\sigma(\zb)\cdot\sigma(\zb)$, from which its invariance
is immediate.  Moreover, it follows from this that the squaring map 
carries the orbits of~$\Spin(10,1)$ to the orbits of~$\SO^\uparrow(10,1)$
and that the image of~$\sigma$ is the union of the origin, the forward
light cone, and the future-directed timelike vectors.
\bigskip

\bigskip
{\bf 6. $\Spin(10,2)$.}  It might be tempting to conjecture that
$\Spin(10,1)$ could be defined directly as the stabilizer of~$p$.  However,
this is not the case, as the stabilizer of~$p$ is larger.
One can see this directly by looking at the alternative expression
$$
p = |\xb_1\w\xb_2|^2 + |\yb_1\w\yb_2|^2     
     - 2\,(\xb_1\cdot\xb_2)\,(\yb_1\cdot\yb_2)
     + 2\,(\xb_1\yb_1)\cdot(\xb_2\yb_2)\ .
$$
which makes it evident that~$p$ is invariant under the 6-dimensional
Lie group
$$
G = \left\{\pmatrix{
   a\,I_8  & 0 &  b\,I_8  &  0\cr
   0 & a'\,I_8  & 0 &  b'\,I_8   \cr 
 c\,I_8    & 0 &  d\,I_8   &   0\cr
 0 &  c'\,I_8 & 0 &     d'\,I_8 }\ 
  \vrule\ ad{-}bc=a'd'{-}b'c'=\pm1\ \right\}\,.
$$
Since $G$ does not lie in~$\Spin(10,1)$, the invariance group of~$p$ must
be properly larger than~$\Spin(10,1)$. 

In particular, consider the $G$-subgroup~$\bR'\simeq\bbR^+$ consisting of 
matrices of the form
$$
\pmatrix{
 t\,I_8  & 0 &  0  &  0\cr
 0 & t^{-1}\,I_8  & 0 &  0  \cr 
 0  & 0 &  t^{-1}\,I_8   &   0\cr
 0 &  0 & 0 &  t\,I_8 }\qquad \hbox{where $t>0$,}
$$
which is not a subgroup of~$\Spin(10,1)$.  Let~$\eur'$ denote its Lie algebra.
Calculation shows that
$$
\eur'\oplus\left[\euspin(10,1),\eur'\right]
= \left\{\pmatrix{
   w\,I_8  & C\,R_\wb &  u\,I_8  &  0\cr
 C\,L_\wb & -w\,I_8  & 0 &  u\,I_8   \cr 
v\,I_8  & 0 &  -w\,I_8 & -C\,R_\wb\cr
 0 &  v\,I_8 & -C\,L_\wb &    w\,I_8   \cr}\ 
  \vrule\ 
\matrix{u,v,w\in\bbR,\cr\noalign{\vskip2pt} \wb\in\bbO}\ \right\}\,.
$$
and that the 
sum~$\euspin(10,1)\oplus\eur'\oplus\left[\euspin(10,1),\eur'\right]$
is closed under Lie bracket.  Thus, this defines a Lie algebra of 
dimension~$66$ that lies in the stabilizer of~$p$.  

The details of further analysis will be omitted, but by using arguments 
similar to those used in previous sections, one sees that this algebra
is isomorphic to~$\euso(10,2)$ and that the connected Lie subgroup of
$\GL_\bbR(\bbO^4)$ whose Lie algebra is this one is simply connected, 
so that it this group is~$\Spin(10,2)$.  Henceforth, this algebra will
be denoted~$\euspin(10,2)$.  Thus,
$$
\euspin(10,2)
= \left\{\pmatrix{
   a_1 +x\,I_8  & C\,R_\wb &  y\,I_8  &  C\,R_\yb\cr
C\,L_\xb & a_3 +w\,I_8  & C\,L_\yb &  u\,I_8   \cr 
z\,I_8    & C\,R_\zb &  a_1 -x\,I_8   &   -C\,R_\xb\cr
 C\,L_\zb &  v\,I_8 & -C\,L_\wb &     a_3 -w\,I_8   \cr}\ 
  \vrule\ 
\matrix{u,v,w,x,y,z\in\bbR,\cr\noalign{\vskip2pt} \wb,\xb,\yb,\zb\in\bbO,\cr
        \noalign{\vskip2pt}\ a\in\euspin(8)}\ \right\}\,.
$$
Moreover, representation theoretic methods show that the only 
connected proper subgroup of~$\SL_\bbR(\bbO^4)$ that properly 
contains~$\Spin(10,2)$ is~$\Sp(16,\bbR)$, the symplectic group 
preserving the symplectic $2$-form~$\Omega$ defined by
$$
\Omega = d\xb_1\,\hat{\cdot}\,d\xb_2 + d\yb_1\,\hat{\cdot}\,d\yb_2\,.
$$
Of course, $\Sp(16,\bbR)$ does not stabilize any nonzero polynomials.
It follows that~$\Spin(10,2)$ is the identity component of the stabilizer 
of~$p$, and hence that the stabilizer of~$p$ must lie in the 
normalizer of~$\Spin(10,2)$ in~$\GL_\bbR(\bbO^4)$.  However, this
normalizer is just~$\bbR^+{\cdot}I_{32}\times\Spin(10,2)$ and the only 
element in~$\bbR^+{\cdot}I_{32}$ that stabilizes~$p$ is the identity
element.  It follows that $\Spin(10,2)$ is the stabilizer of~$p$.

\bigskip
{\bf 7. $\Spin(9,1)$.}  As a final note, inspection reveals that the
subalgebra
$$
\euspin(9,1) = 
\left\{\pmatrix{
   a_1 +x\,I_8  & C\,R_\wb \cr
 C\,L_\xb & a_3 -x\,I_8  }\ 
  \vrule\ 
\matrix{ x\in\bbR,\cr\noalign{\vskip2pt} \wb,\xb\in\bbO,\cr
        \noalign{\vskip2pt}\ a\in\euspin(8)}\ \right\}
\subset\eusl(16,\bbR)\,,
$$
which contains~$\euspin(9)$, is actually the Lie algebra of a faithful
representation of~$\Spin(9,1)$ on~$\bbR^{16}\simeq\bbO^2$.  This action
of~$\Spin(9,1)$ has the interesting feature that it has only two orbits:
The origin and the set of all non-zero vectors.  This follows because the
compact group~$\Spin(9)\subset\Spin(9,1)$ already acts transitively on
the unit spheres, but the larger group does not even preserve the
quadratic form.

\bye